\documentclass[12pt]{article}
\usepackage{amsmath}
\usepackage{mathrsfs}
\usepackage{bbm}
\usepackage{amssymb,a4}
\usepackage{amsfonts,amssymb,amsthm}
\usepackage{extarrows}
\usepackage{color}

\allowdisplaybreaks  % 允许过长的公式分页

\newcommand{\Z}{{\mathbb Z}}
\newcommand{\zd}{{\mathbb Z}^d}

\newcommand{\zp}{{\Z_+}}
\newcommand{\C}{{\mathbb C}}
\newcommand{\cd}{\C^d}

\newcommand{\ad}{\mathrm{ad}}

\newcommand{\cq}{\C_Q}

\newcommand{\dercq}{{\mathrm{Der}(\cq)}}

\newcommand{\zdr}{{\zd\backslash R}}

\newcommand{\g}{\mathfrak{g}}

\newcommand{\derg}{\mathbf{\mathrm{Der}\g}}

\newcommand{\al}{\alpha}
\newcommand{\be}{\beta}
\newcommand{\dt}{\delta}
\newcommand{\lmd}{\lambda}
\newcommand{\gm}{\gamma}
\newcommand{\vf}{\varphi}
\newcommand{\e}{{\epsilon}}
\newcommand{\p}{\partial}
\newcommand{\sgm}{\sigma}

\newcommand{\spanc}[1]{\mathrm{span}_\C\{#1\}}
\newcommand{\sgmf}[2]{(\sgm(#1,#2)-\sgm(#2,#1))}

\newcommand{\bfm}{{\mathbf m}}
\newcommand{\bfn}{{\mathbf n}}
\newcommand{\bfr}{{{\mathbf r}}}
\newcommand{\bfs}{{{\mathbf s}}}

\newcommand{\dergn}{{\derg_\bfn}}

\newcommand{\lm}{L_{\bfm}}
\renewcommand{\ln}{L_{\bfn}}
\newcommand{\lr}{L_{\bfr}}
\newcommand{\ls}{L_{\bfs}}
\newcommand{\lzero}{L_{\bf 0}}

\newcommand{\alm}{\al(\bfm)}
\newcommand{\aln}{\al(\bfn)}
\newcommand{\almn}{\al(\bfm+\bfn)}
\newcommand{\mgm}{{[\bfm]_\gm}}

\newcommand{\tbfr}{{t^{\bfr}}}
\newcommand{\tbfm}{{t^{\bfm}}}
\newcommand{\tbfn}{{t^{\bfn}}}
\newcommand{\tbfs}{{t^{\bfs}}}

\newcommand{\pfend}{\hfill{}$\Box$\end{pf}}

\newtheorem{thm}{Theorem}[section]

\newtheorem{lem}[thm]{Lemma}

\theoremstyle{remark}
\newtheorem*{pf}{{\bf Proof}}

\begin{document}

\begin{center}
{\Large {\bf Automorphisms, derivations and central extensions of
Lie algebras arising from quantum tori}}\\
\vspace{0.5cm}
\end{center}

\begin{center}
{Chengkang Xu\footnote{
The author is supported by the National Natural Science Foundation of China(No. 11626157),
the Science and Technology Foundation of Education Department of
Jiangxi Province(No. GJJ161044).}\\
School of Mathematical Sciences, Shangrao Normal University, Shangrao, Jiangxi,
China}
\end{center}

\begin{abstract}
For a Lie algebra $\g$ related to a quantum torus,
we compute its automorphisms, derivations and universal central extension.
This Lie algebra $\g$ is isomorphic to
a subalgebra of the Lie algebra of derivations over the quantum torus,
and moreover contains a subalgebra isomorphic to the centerless higher rank Virasoro algebra.\\
\noindent
{\bf Keyword}: Quantum torus, automorphism, derivation, central extension,
higher rank Virasoro algebra.
\end{abstract}

\section{Introduction}

\def\theequation{1.\arabic{equation}}
\setcounter{equation}{0}

Throughout this paper, $\C,\Z,\zp$ refer to the set of complex numbers,
integers, and positive integers respectively.
Let $d>1$ be a positive integer,
$\gm=(\gm_1,\gm_2\cdots,\gm_d)$ a $d$-dimensional complex vector,
$Q=(q_{ij})$ a $d\times d$ complex matrix with $q_{ii}=1$
and $q_{ij}q_{ji}=1$ for all $1\leq i, j\leq d$.
Let $\e_1,\e_2\cdots,\e_d$ be a standard $\Z$-basis of $\zd$.
Define a map $\sgm$ on $\zd\times\zd$ by
$$\sgm(\bfm,\bfn)=\prod_{1\leq i<j\leq d}q_{ji}^{m_jn_i},$$
where $\bfm=\sum_{i=1}^d m_i\e_i,\bfn=\sum_{i=1}^d n_i\e_i\in\zd$.
Set
$$R=\{\bfm\in\zd\mid \sgm(\bfm,\bfn)=\sgm(\bfn,\bfm)\text{ for any }\bfn\in\zd\}.$$
The Lie algebra $\g(\gm,Q)$ we consider in this paper is spanned by $\{\lm\mid \bfm\in\zd\}$,
with Lie brackets
$$\begin{aligned}
  &[\lm,\ln]=\sgm(\bfm,\bfn)(\gm\mid\bfn-\bfm)L_{\bfm+\bfn};\\
  &[\lm,\ls]=\sgm(\bfm,\bfs)(\gm\mid\bfs)L_{\bfm+\bfs};\\
  &[\lr,\ls]=(\sgm(\bfr,\bfs)-\sgm(\bfs,\bfr))L_{\bfr+\bfs},
\end{aligned}$$
where $\bfm,\bfn\in R$, $\bfr,\bfs\in\zdr$ and $(\cdot\mid\cdot)$
denotes the inner product on the space $\cd$.

The algebra $\g(\gm,Q)$ is interesting for two reasons.
Firstly it is isomorphic to a subalgebra of the Lie algebra $\dercq$
of derivations of the quantum torus
$\cq=\C[t_1^{\pm1},t_2^{\pm1},\cdots,t_d^{\pm1}]$ related to the matrix $Q$.
In detail, $\cq$ satisfies the commuting relations
$$t_it_i^{-1}=1,\ t_it_j=q_{ij}t_jt_i,\ \ \ \ \ \text{ for any }1\leq i,j\leq d.$$
We denote by $\tbfm$ the monomial $t_1^{m_1}t_2^{m_2}\cdots t_d^{m_d}$,
and by $\p_1,\p_2,\cdots,\p_d$ the degree derivations of $\cq$
corresponding to $t_1,t_2,\cdots, t_d$ respectively.
Then by \cite{BGK} the algebra $\dercq$ has a basis
$$\{\tbfm\p_i,\ad\tbfn\mid \bfm\in R,\bfn\in\zdr, 1\leq i\leq d\},$$
and Lie brackets
$$\begin{aligned}
  &[\tbfm\p_i,\tbfn\p_j]=\sgm(\bfm,\bfn)t^{\bfm+\bfn}(n_i\p_j-m_j\p_i);\\
  &[\tbfm\p_i,\ad\tbfs]=\sgm(\bfm,\bfs)s_it^{\bfm+\bfs};\\
  &[\ad\tbfr,\ad\tbfs]=(\sgm(\bfr,\bfs)-\sgm(\bfs,\bfr))t^{\bfr+\bfs},
\end{aligned}$$
where $\bfm,\bfn\in R$, $\bfr,\bfs\in\zdr$.
One can check that the map defined by
\begin{align*}
  \lm &\mapsto\sum_{i=1}^d\gm_i\tbfm\p_i \text{ if } \bfm\in R;\\
  \ln &\mapsto\ad\tbfn \text{ if } \bfn\in\zdr,
\end{align*}
gives rise to an embedding of the Lie algebra $\g(\gm,Q)$ to $\dercq$.

Secondly, the subalgebra $L=\spanc{\lm\mid\bfm\in R}$ of $\g(\gm,Q)$
is a centerless higher rank Virasoro algebra
if all entries of $Q$ are roots of unity,
in which case the matrix $Q$ and the quantum torus $\cq$ are both called rational.
Denote by $M$ a rank $d$ subgroup of the additive group $\C$.
The higher rank Virasoro algebra $Vir[M]$,
or generalized Virasoro algebra, first introduced in \cite{PZ},
is spanned by $\{e_a, c\mid a\in M\}$ with Lie brackets
$$[c,e_a]=0,\ [e_a,e_b]=(b-a)e_{a+b}+\dt_{a+b,0}\frac{a^3-a}{12}c,\ \ \forall a,b\in M.$$
If $\cq$ is rational, then by \cite{N} and \cite{LZ1}
we know that the subgroup $R$ of $\zd$ is of rank $d$,
and has a $\Z$-basis $k_1\e_1,k_2\e_2,\cdots,k_d\e_d,$
where $1<k_{2i-1}=k_{2i}\in\zp$, $1\leq i\leq z$,
for some positive integer $z$ such that $2z\leq d$,
and $k_{2z+1}=\cdots=k_{d}=1$.
Let $B=diag\{k_1,k_2,\cdots,k_d\}$ and then $R=\{B\bfn\mid\bfn\in\zd\}$.
We call the complex vector $\gm$ generic
if $\gm_1,\gm_2\cdots,\gm_d$ is linearly independent over the field of rational numbers.
If $\gm$ is generic then $M=\sum_{i=1}^d\Z k_i\gm_i$ is a rank $d$ subgroup of $\C$.
Then the map
$$L_{B\bfn}\mapsto e_{(B\gm\mid\bfn)},\ \ \ \bfn\in\zd,$$
defines a Lie algebra isomorphism from $L$ into $Vir[M]$,
whose image is the centerless higher rank Virasoro algebra.
Hence we may think the Lie algebra $\g(\gm,Q)$ as
a "quantum" generalization of the centerless higher rank Virasoro algebra.

We compute the automorphisms,
derivations and universal central extension of $\g(\gm,Q)$ in this paper.
Such computations have always been key problems in the structure theory of Lie algebras.
Many results are known to us,
such as generalized Witt algebras\cite{DZ}, finitely generated graded Lie algebras\cite{Far},
the Lie algebra of skew-derivations over quantum torus\cite{LT},
the twisted Heisenberg-Virasoro algebra\cite{SJ}, and so on.

In this paper we only consider Lie algebra $\g(\gm,Q)$ for the case when $\gm$ is generic.
We remark that $(\gm\mid\bfm)= 0$ if and only if $\bfm={\bf 0}$.
To simplify the notation we denote $\g=\g(\gm,Q)$ if not confused.
The paper is arranged as follows.
Section 2 and Section 3 are devoted to
the computation of automorphisms and derivations of $\g$ respectively.
In the last section we consider central extensions of $\g$ with $Q$ being rational.

\section{Automorphisms}
\def\theequation{2.\arabic{equation}}
\setcounter{equation}{0}

In this section we computer the automorphisms of the algebra $\g$ for generic $\gm$.
Denote
$$\delta(\bfn,R)=\begin{cases}
                      1,&\ \ \text{ if } \bfn\in R;\\
                      0,&\ \ \text{ if } \bfn\notin R.
                 \end{cases}$$
Let $\chi:\zd\longrightarrow\C^*$ be a character of the additive group $\zd$.
Define a linear map $\theta$ on $\g$ by
\begin{equation}\label{eq2.1}
\theta:\ln\mapsto\lmd^{\delta(\bfn,R)}\chi(\bfn)L_{\lmd\bfn},\ \ \ \lmd=\pm1.
\end{equation}
It is easy to check that $\theta$ is a Lie algebra automorphism of $\g$.
The main result in this section is the following
\begin{thm}\label{thm2.1}
 Any Lie algebra automorphism of $\g$ is of the form in (\ref{eq2.1}).
\end{thm}

An element $x\in\g$ is called locally finite if the space
$\spanc{(\ad x)^ny\mid n\in\zp}$
for any $y\in\g$ is finite dimensional.
Clearly the set of locally finite elements in $\g$ is $\C\lzero$.

Let $\theta$ denote an automorphism of $\g$.
Since the image of a locally finite element is still locally finite,
we may assume that
$$\theta(\lzero)=\lmd^{-1}\lzero,\text{ for some nonzero }\lmd\in\C.$$
Moreover, since $\theta$ maps eigenspaces of $\g$ under adjoint action of $\lzero$ into eigenspaces,
we may assume that for any $\bfn\in\zd$,
$$\theta(\ln)=a_\bfn L_{\be(\bfn)},\ \ \ \ a_\bfn\neq 0,\ \ a_{\bf0}=\lmd^{-1},$$
where $\be:\zd\longrightarrow\zd$ is a bijection such that $\be(\bf 0)=\bf 0$.

Applying $\theta$ to $[\lzero,\ln]$ we get
\begin{center}
$\be(\bfn)=\lmd\bfn$ for any $\bfn\in\zd$.
\end{center}
Since $\lmd\bfn\in\zd$, it forces $\lmd=\pm 1$.

\begin{lem}\label{lem2.2}
If $\lmd=1$ then $a_{\bfm+\bfn}=a_\bfm a_\bfn$ for any $\bfm,\bfn\in\zd$.
\end{lem}
\begin{pf}
Note that $\be(\bfn)=\bfn$ for any $\bfn\in\zd$.
Let $\bfm,\bfn\in R$.
Applying $\theta$ to $[\lm,\ln]$ we get
$\sgm(\bfm,\bfn)a_\bfm a_\bfn(\gm\mid\bfn-\bfm)L_{\bfm+\bfn}
 =\sgm(\bfm,\bfn)a_{\bfm+\bfn}(\gm\mid\bfn-\bfm)L_{\bfm+\bfn}$.
So
\begin{equation}\label{eq2.2}
 a_{\bfm+\bfn}=a_\bfm a_\bfn\text{ if } \bfm,\bfn\in R, \bfm\neq\bfn.
\end{equation}
Let $\bfm=\bf0$ in (\ref{eq2.2}) we get $a_{\bf 0}=1=\lmd$.
Let $\bfm+\bfn=\bf0$ in (\ref{eq2.2}) we see
$$ a_{-\bfn}=a_\bfn^{-1}\text{ for any } \bfn\in R.$$
For any $\bfm\in R$, choose $\bfn\in R$ such that $\bfn\neq\bf0,\pm\bfm$.
Then by (\ref{eq2.2}) we have
$$a_{2\bfm}=a_{\bfm-\bfn}a_{\bfm+\bfn}=a_\bfm a_{-\bfn}a_{\bfm}a_\bfn=a_\bfm^2.$$
This proves, together with (\ref{eq2.2}), that
$$a_{\bfm+\bfn}=a_\bfm a_\bfn\text{ if } \bfm,\bfn\in R.$$

Let $\bfm\in R,\bfs\in\zdr$.
Applying $\theta$ to $[\lm,\ls]$ we get
\begin{equation}\label{eq2.3}
 a_{\bfm+\bfs}=a_\bfm a_\bfs\text{ if } \bfm\in R, \bfs\in\zdr.
\end{equation}

For $\bfr, \bfs\in\zdr$ we claim that
$$a_{\bfr+\bfs}=a_\bfr a_\bfs.$$
Apply $\theta$ to $[\lr,\ls]$, and we get
\begin{equation}\label{eq2.4}
(\sgm(\bfr,\bfs)-\sgm(\bfs,\bfr))(a_{\bfr+\bfs}-a_\bfr a_\bfs)=0.
\end{equation}
For $\bfr\in\zdr$, set $G_\bfr=\{\bfn\in\zd\mid \sgm(\bfr,\bfn)=\sgm(\bfn,\bfr)\}$.
Clearly $G_\bfr$ is a proper subgroup of $\zd$ containing $R$.
Choose $\bfn\in\zd\backslash G_\bfr$.
Notice that $\sgm(-\bfr,\bfr+\bfn)\neq\sgm(\bfr+\bfn,-\bfr)$.
So by (\ref{eq2.4}) we have
$a_{\bfr+\bfn}=a_\bfr a_\bfn=a_\bfr a_{-\bfr}a_{\bfr+\bfn}$.
Therefore
\begin{equation}\label{eq2.5}
a_{-\bfr}a_\bfr=1=a_{\bf0},\ \ a_{-\bfr}=a_\bfr^{-1}\text{ for any } \bfr\in\zdr.
\end{equation}
Now we may prove the claim in three cases.\\
{\bf Case 1}: $\sgm(\bfr,\bfs)\neq\sgm(\bfs,\bfr)$.
The claim is obvious by (\ref{eq2.4}).\\
{\bf Case 2}: $\bfr+\bfs\in R\backslash\{\bf 0\}$(this implies $\sgm(\bfr,\bfs)=\sgm(\bfs,\bfr)$).
By (\ref{eq2.3}) and (\ref{eq2.5}) we have
$$a_\bfr=a_{\bfr+\bfs-\bfs}=a_{\bfr+\bfs}a_{-\bfs}=a_{\bfr+\bfs}a_{\bfs}^{-1},$$
and the claim follows in this case.\\
{\bf Case 3}: $\sgm(\bfr,\bfs)=\sgm(\bfs,\bfr)$ and $\bfr+\bfs\notin R$.
A basic group theory shows that $G_{\bfr+\bfs}\cup G_\bfr \cup G_\bfs\neq \zd$.
We may choose $\bfn\in\zd$ and $\bfn\notin G_{\bfr+\bfs}\cup G_\bfr \cup G_\bfs$.
Noticing that
$$\sgm(\bfr+\bfs+\bfn,-\bfn)\neq\sgm(-\bfn,\bfr+\bfs+\bfn),
\sgm(\bfr,\bfs+\bfn)\neq\sgm(\bfs+\bfn,\bfr) \text{ and }
\sgm(\bfs+\bfn,-\bfn)\neq\sgm(-\bfn,\bfs+\bfn),$$
we have, by Case 1,
$$
a_{\bfr+\bfs}=a_{\bfr+\bfs+\bfn}a_{-\bfn}=a_\bfr a_{\bfs+\bfn}a_{-\bfn}=a_\bfr a_{\bfs}.
$$
This proves the claim, and hence the lemma.
\pfend

For $\bfm,\bfn\in\zd$, denote by $\iota(\bfm,\bfn)$ the number of elements in
$\{\bfm,\bfn,\bfm+\bfn\}\cap R$.

\begin{lem}\label{lem2.3}
If $\lmd=-1$ then for any $\bfm,\bfn\in\zd$,
$$\begin{cases}
a_{\bfm+\bfn}=a_\bfm a_\bfn,\ \  &\text{if }\iota(\bfm,\bfn)=0;\\
a_{\bfm+\bfn}=-a_\bfm a_\bfn,\ \ &\text{ if }\iota(\bfm,\bfn)>0.
\end{cases}$$
\end{lem}
\begin{pf}
The proof is parallel to that of Lemma \ref{lem2.2}.
\pfend

Define $\chi:\bfn\mapsto\lmd^{\delta(\bfn,R)}a_\bfn$.
By Lemma \ref{lem2.2} and Lemma \ref{lem2.3}
one can check that $\chi $ is a character of $\zd$.
So
$$ a_\bfn=\lmd^{\delta(\bfn,R)}\chi(\bfn),$$
which proves Theorem \ref{thm2.1}.

\section{Derivations}

\def\theequation{3.\arabic{equation}}
\setcounter{equation}{0}

A linear map $D:\g\longrightarrow\g$ is called a derivation of $\g$ if
$$D[\lm,\ln]=[D\lm,\ln]+[\lm,D\ln]\text{ for any }\bfm,\bfn\in\zd.$$
In this section we computer derivations of $\g$ for generic $\gm$.
Since $\g$ is $\zd$-graded and finitely generated,
a result from \cite{Far} shows that the Lie algebra $\derg$
of derivations of $\g$ is also $\zd$-graded.
Specifically we have the following
\begin{thm}
The $\zd$-graded structure of $\derg$ is $\derg=\bigoplus_{\bfn\in\zd}\dergn$ where
$$\dergn=\begin{cases}
            \spanc{\p_1,\p_2,\cdots,\p_d}, &\text{ if }\bfn=\bf0;\\
            \spanc{\ad\ln}, &\text{ if }\bfn\neq\bf0,
         \end{cases}$$
where $\p_i:\lm\mapsto m_i\lm$, $1\leq i\leq d$ for $\bfm=(m_1,\cdots,m_d)\in\zd$.
\end{thm}
\begin{pf}
Since $\p_i, 1\leq i\leq d$, are derivations of degree $\bf0$ on $\g$,
and linearly independent,
to prove $\derg_{\bf0}=\spanc{\p_1,\p_2,\cdots,\p_d}$
it suffices to show that $\dim\derg_{\bf0}\leq d$.
Let $D\in\derg_{\bf0}$ and
$$D(\lm)=\vf(\bfm)\lm,\ \ \ \bfm\in\zd$$
for some function $\vf:\zd\longrightarrow\C$.

Let $\bfm,\bfn\in R$. Apply $D$ to $[\lm,\ln]$ and we get
$$(\gm\mid\bfn-\bfm)(\vf(\bfm+\bfn)-\vf(\bfm)-\vf(\bfn))=0,$$
which shows
\begin{equation}\label{eq3.1}
\vf(\bfm+\bfn)=\vf(\bfm)+\vf(\bfn)\text{ for } \bfm\neq\bfn\in R.
\end{equation}
Let $\bfn=\bf0$ in (\ref{eq3.1}) we see $\vf({\bf0})= 0$.
Now for $\bfm\neq\bf0$, choose $\bfs\in R$ such that $\bfs\neq {\bf0}, \bfm, \frac{1}{2}\bfm$,
then by (\ref{eq3.1}) we have
$$\vf(2\bfm)=\vf(2\bfm-\bfs)+\vf(\bfs)=\vf(\bfm)+\vf(\bfm-\bfs)+\vf(\bfs)=2\vf(\bfm).$$
This, together with (\ref{eq3.1}), shows
\begin{equation}\label{eq3.2}
\vf(\bfm+\bfn)=\vf(\bfm)+\vf(\bfn)\text{ for any } \bfm,\bfn\in R.
\end{equation}

Applying $D$ to $[\lm,\ls]$ for $\bfm\in R, \bfs\in\zdr$, we get
\begin{equation}\label{eq3.3}
\vf(\bfm+\bfs)=\vf(\bfm)+\vf(\bfs)\text{ for any } \bfm\in R, \bfs\in\zdr.
\end{equation}

For $\bfr,\bfs\in\zdr$ we claim that $\vf(\bfr+\bfs)=\vf(\bfr)+\vf(\bfs)$.
From $D[\lr,\ls]=[D\lr,\ls]+[\lr,D\ls]$ we get
$$(\sgm(\bfr,\bfs)-\sgm(\bfs,\bfr))(\vf(\bfr+\bfs)-\vf(\bfr)-\vf(\bfs))=0,$$
which implies
\begin{equation}\label{eq3.4}
\vf(\bfr+\bfs)=\vf(\bfr)+\vf(\bfs)\text{ if } \sgm(\bfr,\bfs)\neq\sgm(\bfs,\bfr).
\end{equation}
Let $\bfn\in\zd\backslash G_\bfr$.
Noticing that $\sgm(-\bfr,\bfr+\bfn)\neq\sgm(\bfr+\bfn,-\bfr)$, we have
$$\vf(\bfr+\bfn)=\vf(\bfr)+\vf(\bfn)=\vf(\bfr)+\vf(-\bfr)+\vf(\bfr+\bfn).$$
Therefore
\begin{equation}\label{eq3.5}
\vf(-\bfr)+\vf(\bfr)=0=\vf({\bf0}),\ \ \vf(-\bfr)=-\vf(\bfr)\text{ for any } \bfr\in\zdr.
\end{equation}
If $\bfr+\bfs\in R\backslash\{\bf0\}$(this implies $\sgm(\bfr,\bfs)=\sgm(\bfs,\bfr)$),
then by (\ref{eq3.3}) and (\ref{eq3.5}) we have
\begin{equation}\label{eq3.6}
\vf(\bfr+\bfs)=\vf(\bfr+\bfs)+\vf(-\bfs)+\vf(\bfs)=\vf(\bfr)+\vf(\bfs).
\end{equation}
If $\bfr+\bfs\in\zdr$ and $\sgm(\bfr,\bfs)=\sgm(\bfs,\bfr)$,
then we choose $\bfn\in\zd$ such that $\bfn\notin G_\bfr\cup G_\bfs\cup G_{\bfr+\bfs}$,
and we have
$$\vf(\bfr+\bfs)=\vf(\bfr+\bfs+\bfn)+\vf(-\bfn)
=\vf(\bfr)+\vf(\bfn+\bfs)+\vf(-\bfn)=\vf(\bfr)+\vf(\bfs),$$
which, together with (\ref{eq3.4}),(\ref{eq3.5}) and (\ref{eq3.6}), proves the claim.
Hence
$$\vf(\bfm+\bfn)=\vf(\bfm)+\vf(\bfn)\text{ for any } \bfm,\bfn\in\zd,$$
and $\vf(\bfm)=\sum_{i=1}^d m_i\vf(\e_i).$
So $\dim\derg_{\bf0}\leq d$ and $\derg_{\bf0}=\spanc{\p_1,\p_2,\cdots,\p_d}$.

For $\bfn\neq{\bf0}$, suppose $D\in\dergn$ and
write $D(\lm)=\phi(\bfm)L_{\bfm+\bfn}$ for any $\bfm\in\zd$,
where $\phi:\zd\longrightarrow\C$.
Apply $D$ to $[\lzero,\lm]$, we get
$$\phi({\bf0})[\ln,\lm]+\phi(\bfm)(\gm\mid \bfm+\bfn)L_{\bfm+\bfn}
   =\phi(\bfm)(\gm\mid \bfm)L_{\bfm+\bfn},$$
which implies
$$\phi(\bfm)=\phi({\bf0})\frac{(\gm\mid\bfm-\bfn)}{(\gm\mid \bfn)},$$
which implies that $\dim\dergn\leq1$.
Moreover, $\dergn$ contains the inner derivation $\ad\ln$.
So $\dergn=\spanc{\ad\ln}$.
\pfend

At last we remark that the inner derivation
$\ad\lzero=\gm_1\p_1+\gm_2\p_2+\cdots+\gm_d\p_d$.

\section{Universal central extension}

\def\theequation{4.\arabic{equation}}
\setcounter{equation}{0}

Since $[\g,\g]=\g$, the Lie algebra $\g$ has a universal central extension.
In this section we compute the universal central extension of $\g$
for rational matrix $Q$ and generic $\gm$.
As mentioned in Section 1 we consider $\g$ as a subalgebra of the Lie algebra $\dercq$,
and we need some basics about the rational quantum torus $\cq$ and the map $\sgm$.
The following lemma is from \cite{N} and \cite{LZ1}.
\begin{lem}\label{lem4.1}
Up to an isomorphism of $\cq$, we may assume that all entries of the matrix $Q$ are 1
except $q_{2i-1}=q_{2i-1,2i}$ and $q_{2i}=q_{2i,2i-1}={q_{2i-1}}^{-1}$ for $1\leq i\leq z$,
where $z\in\zp$ with $2z\leq d$,
and the order $k_i$ of $q_i$, $1\leq i\leq 2z$,
as roots of unity, satisfy $k_{i+1}\mid k_i$ and
$k_{2j-1}=k_{2j}$ for $1\leq i\leq 2z-1,\ 1\leq j\leq z$. Let $k_l=1$ for $l>2z$.
Then the subgroup $R$ has the form
$$R=\bigoplus_{i=1}^d k_i\Z\e_{i}. $$
\end{lem}

From now on we will always assume that
the numbers $k_1,\cdots,k_d$ are fixed, $Q$ and $R$ have the form as in Lemma \ref{lem4.1}.

\begin{lem}\label{lem4.2}
(1) $\sgm(\bfm+\bfn,\bfr+\bfs)=\sgm(\bfm,\bfr)\sgm(\bfm,\bfs)\sgm(\bfn,\bfr)\sgm(\bfn,\bfs)$.\\
(2) $\sgm(\bfm,\bfn)\sgm(\bfm,-\bfn)=1$.\\
(3) $\sgm(-\bfm,-\bfn)=\sgm(\bfm,\bfn)$.\\
(4) $\sgm(\bfm,\bfn)=1$ for all $\bfm\in R,\bfn\in\zd$.\\
(5) $\sgm(\e_i,\e_j)=\begin{cases}
               q_{i,i-1}  &\text{ if }i\text{ is even and }j=i-1;\\
               1          &\text{ otherwise}.
               \end{cases}$

\end{lem}

The lemma is easy to be checked
and critical to the following computation.
The Lie bracket of $\g$ becomes
$$\begin{aligned}
  &[\lm,\ln]=(\gm\mid\bfn-\bfm)L_{\bfm+\bfn};\\
  &[\lm,\ls]=(\gm\mid\bfs)L_{\bfm+\bfs};\\
  &[\lr,\ls]=(\sgm(\bfr,\bfs)-\sgm(\bfs,\bfr))L_{\bfr+\bfs},
\end{aligned}$$

Now we start to compute the central extensions of $\g$.
Let $\al:\g\times\g\longrightarrow\C$ be an arbitrary 2-cocycle of $\g$.
Hence
\begin{equation}\label{eq4.1}
 \al([\lm,\ln],\ls)+\al([\ls,\lm],\ln)+\al([\ln,\ls],\lm)=0,
 \text{ for any }\bfm,\bfn,\bfs\in\zd.
\end{equation}
Define a linear function $f_\al$ on $\g$ by
$$\begin{cases}
 f_\al(\lm)=\frac{1}{(\gm\mid2k_1\e_1-\bfm)}\al(L_{\bfm-k_1\e_1},L_{k_1\e_1})
                        &\text{ for }\bfm\in R,\bfm\neq 2k_1\e_1;\\
 f_\al(L_{2k_1\e_1})=\frac{1}{(\gm\mid2k_1\e_1)}\al(\lzero,L_{2k_1\e_1});& \\
 f_\al(\ls)=\frac{1}{(\gm\mid\bfs)}\al(\lzero,\ls)&\text{ for }\bfs\in\zdr,
\end{cases}$$
and a 2-coboundary $\psi_{f_\al}$ by
$$\psi_{f_\al}(\lm,\ln)=f_\al([\lm,\ln]).$$
It is easy to check that
\begin{align*}
 &\psi_{f_\al}(L_{\bfm-k_1\e_1},L_{k_1\e_1})=\al(L_{\bfm-k_1\e_1},L_{k_1\e_1})
    \text{ for }\bfm\in R,\bfm\neq 2k_1\e_1\\
 &\psi_{f_\al}(\lzero,L_{2k_1\e_1})=\al(\lzero,L_{2k_1\e_1}),\\
 &\psi_{f_\al}(\lzero,\ls)=\al(\lzero,\ls)\text{ for }\bfs\in\zdr.
\end{align*}
Since $\al-\psi_{f_\al}$ is equivalent to $\al$, we may assume
$$\begin{cases}
 \al(L_{\bfm-k_1\e_1},L_{k_1\e_1})=0&\text{ for }\bfm\in R,\bfm\neq 2k_1\e_1;\\
 \al(\lzero,L_{2k_1\e_1})=0;& \\
 \al(\lzero,\ls)=0&\text{ for }\bfs\in\zdr.
\end{cases}$$

\begin{lem}\label{lem4.3}
$\al(\lm,\ln)=0$ if $\bfm+\bfn\neq{\bf0}$.
\end{lem}
\begin{pf}
Take $\bfm,\bfn,\bfs\in R$ in (\ref{eq4.1}) and we get
\begin{equation}\label{eq4.2}
 (\gm\mid\bfn-\bfm)\al(L_{\bfm+\bfn},\ls)+(\gm\mid\bfs-\bfn)\al(L_{\bfs+\bfn},\lm)
   +(\gm\mid\bfm-\bfs)\al(L_{\bfm+\bfs},\ln)=0.
\end{equation}
Let $\bfs={\bf0}$ in (\ref{eq4.2}) and we have
\begin{equation}\label{eq4.3}
  (\gm\mid\bfn+\bfm)\al(\lm,\ln)=(\gm\mid\bfm-\bfn)\al(L_{\bfm+\bfn},\lzero).
\end{equation}
Replacing $\bfm$ by $\bfm-k_1\e_1$ and $\bfn$ by $k_1\e_1$, we get
\begin{equation}\label{eq4.4}
  (\gm\mid\bfm-2k_1\e_1)\al(\lm,\lzero)=(\gm\mid\bfm)\al(L_{\bfm-k_1\e_1},L_{k_1\e_1})=0,
\end{equation}
which implies
\begin{equation}\label{eq4.5}
  \al(\lm,\lzero)=0 \text{ for any }\bfm\in R.
\end{equation}
So by (\ref{eq4.3}) and (\ref{eq4.5}), the lemma stands for the case when $\bfm,\bfn\in R$.

Let $\bfm\in R,\bfn={\bf0},\bfs\in\zdr$ in (\ref{eq4.1}), we get
$$(\gm\mid\bfm+\bfs)\al(\lm,\ls)=-(\gm\mid\bfs)\al(L_{\bfm+\bfs},L_{\bf0})=0,$$
which shows
$$\al(\lm,\ls)=0\text{ for any }\bfm\in R,\bfs\in\zdr.$$

For $\bfr,\bfs\in\zdr$, since
$$(\gm\mid\bfs)\al(\lr,\ls)=\al(\lr,[\lzero,\ls])
   =\al([\lr,\lzero],\ls)+\al(\lzero,[\lr,\ls])=-(\gm\mid\bfr)\al(\lr,\ls),$$
we get
$$(\gm\mid\bfr+\bfs)\al(\lr,\ls)=0,$$
which proves the lemma for the case when $\bfr,\bfs\in\zdr$.
\pfend

Now we only need to determine $\al(\lm,\ln)$ when $\bfm+\bfn=\bf0$.
To simplify the notation,
from now on we denote $\alm=\al(\lm,L_{-\bfm})$ for any $\bfm\in\zd$.

\begin{lem}\label{lem4.4}
Let $\bfm\in R$, then $\alm=\frac{\mgm^3-\mgm}{6}\al(2k_1\e_1)$,
where $\mgm=\frac{(\gm\mid\bfm)}{(\gm\mid k_1\e_1)}$.
\end{lem}
\begin{pf}
Let $\bfm,\bfn\in R$ and $s=-\bfm-\bfn$ in (\ref{eq4.1}), we get
\begin{equation}\label{eq4.6}
(\gm\mid\bfn-\bfm)\al(\bfm+\bfn)=(\gm\mid 2\bfm+\bfn)\aln-(\gm\mid 2\bfn+\bfm)\alm
\text{ for any }\bfm,\bfn\in R.
\end{equation}
Let $\bfn=-k_1\e_1$ in (\ref{eq4.6}) and notice that $\al(k_1\e_1)=0$, we have
\begin{equation}\label{eq4.7}
(\gm\mid\bfm-2k_1\e_1)\alm=(\gm\mid\bfm+k_1\e_1)\al(\bfm-k_1\e_1)\text{ for any }\bfm\in R.
\text{ for any }\bfm,\bfn\in R.
\end{equation}

First we prove the lemma for the case $\bfm=l k_1\e_1,\ l\in\Z$.
In this case $\mgm=l$.
If $l>2$, by (\ref{eq4.7}), we have
\begin{align*}
\al(lk_1\e_1)&=\frac{(\gm\mid lk_1\e_1+k_1\e_1)}{(\gm\mid lk_1\e_1-2k_1\e_1)}
              \frac{(\gm\mid lk_1\e_1)}{(\gm\mid lk_1\e_1-3k_1\e_1)}
              \frac{(\gm\mid lk_1\e_1-k_1\e_1)}{(\gm\mid lk_1\e_1-4k_1\e_1)}\cdots
              \frac{(\gm\mid 4k_1\e_1)}{(\gm\mid k_1\e_1)}\al(2k_1\e_1)\\
             &=\frac{l^3-l}{6}\al(2k_1\e_1).
\end{align*}
If $l<-2$, we have
$$\al(lk_1\e_1)=-\al(-lk_1\e_1)=-\frac{(-l)^3-(-l)}{6}\al(2k_1\e_1)=\frac{l^3-l}{6}\al(2k_1\e_1).$$
If $-2\leq l\leq 2$, the expression can be easily checked for each $l$.

Secondly, we deal with the case $\bfm=l\e_i\in R,\ i>1$.
Set $\bfm=2k_1\e_1,\bfn=l\e_i$ in (\ref{eq4.6})
and by (\ref{eq4.7}) we have
\begin{align*}
  &(\gm\mid 4k_1\e_1+l\e_i)\al(l\e_i)-(\gm\mid 2k_1\e_1+2l\e_i)\al(2k_1\e_1)\\
 =&(\gm\mid l\e_i-2k_1\e_1)\al(l\e_i+2k_1\e_1)\\
 =&(\gm\mid l\e_i-2k_1\e_1)\frac{(\gm\mid l\e_i+3k_1\e_1)}{(\gm\mid l\e_i)}
    \frac{(\gm\mid l\e_i+2k_1\e_1)}{(\gm\mid l\e_i-k_1\e_1)}\al(l\e_i),
\end{align*}
which gives
\begin{equation}\label{eq4.8}
  \al(l\e_i)=\frac{x^3-x}{6}\al(2k_1\e_1),\ \
  x=\frac{(\gm\mid l\e_i)}{(\gm\mid k_1\e_1)}=[l\e_i]_{\gm},
\end{equation}
proving the lemma for the case $\bfm=l\e_i\in R,\ i>1$.

At last we treat the general case.
Let $k$ denote the maximal index in $\bfm=(m_1,\cdots,m_d)$ such that $m_k\neq0$.
We use induction on $k$.
Denote $x_i=\frac{(\gm\mid m_i\e_i)}{(\gm\mid k_1\e_1)}$ for $1\leq i\leq k$,
and we have $\sum\limits_{i=1}^kx_i=\mgm$.
The case $k=1$ is just the special case $\bfm\in\Z k_1e_1$.
Suppose the lemma stands for any $\bfn\in R$ such that $n_k=n_{k+1}=\cdots=n_d=0$.
Replacing $\bfm$ by $\bfm-m_k\e_k$ and $\bfn$ by $m_k\e_k$ in (\ref{eq4.6}) gives
$$
 (\gm\mid2m_k\e_k-\bfm)\al(\bfm)=(\gm\mid 2\bfm-m_k\e_k)\al(m_k\e_k)
       -(\gm\mid m_k\e_k+\bfm)\al(\bfm-m_k\e_k).
$$
Hence by (\ref{eq4.8}) and the inductional hypothesis we get
\begin{align*}
\alm&=-\frac{2\mgm-x_k}{\mgm-2x_k}\al(m_k\e_k)
      +\frac{2\mgm+x_k}{\mgm-2x_k}\al(\bfm-m_k\e_k)\\
    &=-\frac{2\mgm-x_k}{\mgm-2x_k}\frac{x_k^3-x_k}{6}\al(2k_1\e_1)
      +\frac{2\mgm+x_k}{\mgm-2x_k}\left((\sum_{j=1}^{k-1}x_j)^3-\sum_{j=1}^{k-1}x_j\right)
       \frac{1}{6}\al(2k_1\e_1)\\
    &=\frac{1}{6}\left((\sum_{j=1}^{k}x_j)^3-\sum_{j=1}^{k}x_j\right)\al(2k_1\e_1)\\
    &=\frac{\mgm^3-\mgm}{6}\al(2k_1\e_1),
\end{align*}
proving the lemma.
\pfend

Let $\bfm\in R$, $\bfn\in\zdr$ and $\bfs=-\bfm-\bfn$ in (\ref{eq4.1}), we have
\begin{equation}\label{eq4.9}
(\gm\mid\bfn)\al(\bfm+\bfn)=(\gm\mid\bfm+\bfn)\aln,
\end{equation}
which means, together with Lemma \ref{lem4.3}, that
$\al(\lm,\ln),\bfm,\bfn\in\zdr$ is determined by
$$\alm,\ \ \ \bfm\in\Gamma=\{\bfn\in\zdr\mid 0\leq n_i<k_i,1\leq i\leq d\}.$$

Now let $\bfm,\bfn\in\zdr$ and $\bfs=-\bfm-\bfn$ in (\ref{eq4.1}), we get
\begin{equation}\label{eq4.10}
 \sgmf{\bfm}{\bfn}\almn=\sgmf{\bfn}{-\bfm}\left(\sgm(\bfn,-\bfn)\alm+\sgm(\bfm,-\bfm)\aln\right)
\end{equation}
\begin{lem}\label{lem4.5}
Let $\bfm\in\zdr$. Suppose $m_i\neq0$ for some odd $1\leq i\leq 2z$,
then
$$\alm=\frac{(\gm\mid\bfm)}{(\gm\mid\e_1)}\sgm(\bfm,-\bfm)\al(\e_1).$$
\end{lem}
\begin{pf}
With suitable substitutions in (\ref{eq4.10}) we get three equations as follows.
\begin{align*}
\alm  ={q_{i+1,i}}^{-m_i}\al(\bfm-\e_{i+1})&+\sgm(\bfm,-\bfm)\al(\e_{i+1})
       \ \ \ \text{ if } k_i\nmid m_i;\\
\al(\bfm-\e_{i+1})  ={q_{i+1,i}}^{-k(m_{i+1}-1)}\al(\bfm-\e_{i+1}-k\e_i)
                      &+\sgm(\bfm,-\bfm){q_{i+1,i}}^{m_i}\al(k\e_i)\\
               &\text{ if } k_{i+1}\nmid(m_{i+1}-1)k;\\
\al(\bfm-\e_{i+1}-k\e_i)={q_{i+1,i}}^{m_i-k}\al(\bfm-k\e_i)
                       &-\sgm(\bfm,-\bfm){q_{i+1,i}}^{k(m_{i+1}-1)+m_i}\al(\e_{i+1})\\
              &\text{ if } k_i\nmid (m_i-k).
\end{align*}
From these three equations it follows that
\begin{equation}\label{eq4.11}
\alm={q_{i+1,i}}^{-km_{i+1}}\al(\bfm-k\e_i)+\sgm(\bfm,-\bfm)\al(k\e_i),
\end{equation}
if $k_i\nmid m_i$ and $k$ satisfies $\ k_{i+1}\nmid(m_{i+1}-1)k$ and $k_i\nmid (m_i-k)$.
Particularly we have
$$
\alm={q_{i+1,i}}^{-m_{i+1}}\al(\bfm-\e_i)+\sgm(\bfm,-\bfm)\al(\e_i)
     \text{ for any }\bfm\in\Gamma.
$$
And using induction on $m_i$ we get
\begin{equation}\label{eq4.12}
\alm={q_{i+1,i}}^{-m_{i+1}(m_{i}-1)}\al(\bfm-(m_{i}-1)\e_i)+\sgm(\bfm,-\bfm)(m_{i}-1)\al(\e_i)
  \text{ for any }\bfm\in\Gamma.
\end{equation}
In particular,
\begin{equation}\label{eq4.13}
\al(l\e_i)=l\al(\e_i)\text{  for } 0<l<k_i.
\end{equation}

First we deal with the special case when $\bfm\in\Gamma$
with $m_i\neq0$ for some odd $1\leq i\leq2z$.
Set $\bfn=\bfm-(m_{i}-1)\e_i$ and notice that $n_i=1, k_i\e_i\in R$.
Take $0<l<k_i$ satisfies the conditions in (\ref{eq4.11}).
Notice that $k_i-l$ still satisfies these conditions.
Replacing $\bfm$ by $\bfn+k_i\e_i$ and taking $k=l$ in (\ref{eq4.11}),
then replacing $\bfm$ by $\bfn+(k_i-l)\e_i$ and taking $k=k_i-l$ in (\ref{eq4.11}),
and using (\ref{eq4.13}), we get
\begin{align*}
\al(\bfn+k_i\e_i)=&{q_{i+1,i}}^{-ln_{i+1}}\al(\bfn+(k_i-l)\e_i)+\sgm(\bfn,-\bfn)\al(l\e_i)\\
                 =&{q_{i+1,i}}^{-ln_{i+1}}\left({q_{i+1,i}}^{-n_{i+1}(k_i-l)}\aln
                   +\sgm(\bfn-l\e_i,-\bfn+l\e_i)\al((k_i-l)\e_i)\right)\\
                  &   +\sgm(\bfn,-\bfn)\al(l\e_i)\\
                 =&\aln+k_i\sgm(\bfn,-\bfn)\al(\e_i).
\end{align*}
Using (\ref{eq4.9}) we see
\begin{align*}
(\gm\mid\bfn+k_i\e_i)\aln=(\gm\mid\bfn)\al(\bfn+k_i\e_i)
=(\gm\mid\bfn)\aln+(\gm\mid\bfn)k_i\sgm(\bfn,-\bfn)\al(\e_i).
\end{align*}
Hence
$$\aln=\frac{(\gm\mid\bfn)}{(\gm\mid\e_i)}\sgm(\bfn,-\bfn)\al(\e_i).$$
This, together with (\ref{eq4.12}), implies
\begin{equation}\label{eq4.14}
\alm=\frac{(\gm\mid\bfm)}{(\gm\mid\e_i)}\sgm(\bfm,-\bfm)\al(\e_i).
\end{equation}
Especially take $\bfm=\e_1+\e_i$ and we see
$$\frac{(\gm\mid\bfm)}{(\gm\mid\e_i)}\sgm(\bfm,-\bfm)\al(\e_i)=\alm=
  \frac{(\gm\mid\bfm)}{(\gm\mid\e_1)}\sgm(\bfm,-\bfm)\al(\e_1).$$
Therefore,
$$\al(\e_i)=\frac{(\gm\mid\e_i)}{(\gm\mid\e_1)}\al(\e_1).$$
Substitute this into (\ref{eq4.14})
and we have proved the lemma for the case $\bfm\in\Gamma$ with $m_i\neq0$
for some odd $1\leq i\leq2z$.

Now we prove the lemma for the general case.
Let $\bfr\in R$ such that $\bfm-\bfr\in\Gamma$ satisfies
$m_i- r_i\neq 0$ for some odd $1\leq i\leq 2z$.
Replacing $\bfm$ by $\bfr$, $\bfn$ by $\bfm-\bfr$ in (\ref{eq4.9}),
and applying the special case above, we see
$$
 (\gm\mid\bfm-\bfr)\alm=(\gm\mid\bfm)\al(\bfm-\bfr)
   =(\gm\mid\bfm)\frac{(\gm\mid\bfm-\bfr)}{(\gm\mid\e_1)}
   \sgm(\bfm-\bfr,-\bfm+\bfr)\al(\e_1).
$$
So
$$\alm=\frac{(\gm\mid\bfm)}{(\gm\mid\e_1)}\sgm(\bfm-\bfr,-\bfm+\bfr)\al(\e_1)
      =\frac{(\gm\mid\bfm)}{(\gm\mid\e_1)}\sgm(\bfm,-\bfm)\al(\e_1),$$
proving the lemma.

\pfend

\begin{lem}\label{lem4.6}
Let $\bfm\in\zdr$. Then
$$\alm=\frac{(\gm\mid\bfm)}{(\gm\mid\e_1)}\sgm(\bfm,-\bfm)\al(\e_1).$$
\end{lem}
\begin{pf}
By Lemma \ref{lem4.5} we only need to prove the case $m_{2i-1}=0$ for all $1\leq i\leq z$.
Since $\bfm\in\zdr$, there exists some $1\leq i\leq z$ such that $m_{2i}\neq0$.
Replacing $\bfm$ by $\bfm-\e_{2i-1}$ and $\bfn$ by $\e_{2i-1}$ in (\ref{eq4.10}),
using Lemma \ref{lem4.5}, we get
\begin{align*}
\alm &={q_{2i,2i-1}}^{-m_{2i}}(\al(\bfm-\e_{2i-1})
                      +\sgm(\bfm-\e_{2i-1},-\bfm+\e_{2i-1})\al(\e_{2i-1}))\\
     &={q_{2i,2i-1}}^{-m_{2i}}\sgm(\bfm-\e_{2i-1},-\bfm+\e_{2i-1})\left(
           \frac{(\gm\mid\bfm-\e_{2i-1})}{(\gm\mid \e_1)}\al(\e_1)
                      +\frac{(\gm\mid \e_{2i-1})}{(\gm\mid \e_1)}\al(\e_1)\right)\\
     &=\sgm(\bfm,-\bfm)\frac{(\gm\mid\bfm)}{(\gm\mid\e_1)}\al(\e_1).
\end{align*}
\pfend

Now we may summarise our result in this section,
which follows directly from Lemma \ref{lem4.3}, Lemma \ref{lem4.4} and Lemma \ref{lem4.6}.
\begin{thm}
The universal central extension $\tilde\g$ of the Lie algebra $\g$ has a center of dimension 2,
which has basis $c_1=2\al(2k_1\e_1),\ c_2=k_1\al(\e_1)$.
More specifically, the Lie algebra $\tilde\g$ has Lie bracket
$$\begin{aligned}
  &[\lm,\ln]'=(\gm\mid\bfn-\bfm)L_{\bfm+\bfn}
        +\delta_{\bfm+\bfn,{\bf0}}\frac{\mgm^3-\mgm}{12}c_1;\\
  &[\lm,\ls]'=(\gm\mid\bfs)L_{\bfm+\bfs};\\
  &[\lr,\ls]'=(\sgm(\bfr,\bfs)-\sgm(\bfs,\bfr))L_{\bfr+\bfs}
       +\delta_{\bfr+\bfs,{\bf0}}\sgm(\bfr,-\bfr)[\bfr]_\gm c_2;\\
  &[c_1,\tilde\g]'=[c_2,\tilde\g]'=0.
\end{aligned}$$
where $\bfm,\bfn\in R$, $\bfr,\bfs\in\zdr$,
$\mgm=\frac{(\gm\mid\bfm)}{(\gm\mid k_1\e_1)}$
and $[\bfr]_\gm=\frac{(\gm\mid\bfr)}{(\gm\mid k_1\e_1)}$.
\end{thm}

At last we mention that the subalgebra
$\spanc{\lm,c_1\mid \bfm\in R}$ of $\tilde\g$
is a higher rank Virasoro algebra.

\end{document}